\documentclass{amsart}     
\usepackage{amssymb,epsfig,psfrag,subfigure} 
\usepackage{amsmath,amscd} 
\newtheorem{theorem}{Theorem}[section]    
\newtheorem{lemma}[theorem]{Lemma}         
\newtheorem{corollary}[theorem]{Corollary}

\theoremstyle{definition}
\newcommand{\ZZ}{\mathbb{Z}}
\newcommand{\NN}{\mathbb{N}}
\newcommand{\QQ}{\mathbb{Q}}
\newcommand{\CC}{\mathbb{C}}
\newcommand{\HH}{\mathbb{H}}
\newcommand{\RR}{\mathbb{R}}

\def\a{\alpha}
\def\b{\beta}
\def\g{\gamma}

\def\G{\Gamma}

\def\e{\epsilon}
\def\D{\Delta}

\def\PSL{\mbox{\rm{PSL}}}

\def\d{\mbox{\rm{d}}}

\def\SL{\mbox{\rm{SL}}}

\def\Isom{\mbox{Isom}}

\def\Vol{\mbox{\rm{Vol}}}
\def\Area{\mbox{\rm{Area}}}
\def\tube{\mbox{\rm{tube}}}

\def\d{\partial}

\def\S{\Sigma}

\begin{document}

\title{Bounds on exceptional Dehn filling II}                    
\author{Ian Agol}                  
\address{University of California, Berkeley \\ 970 Evans Hall \#3840 \\ Berkeley, CA 94720-3840}                  


\email{ianagol@math.berkeley.edu}      
\thanks{partially supported by NSF grant DMS-0504975}               
%

\begin{abstract}   
We show that there are at most finitely many one cusped orientable hyperbolic
3-manifolds which have more than eight non-hyperbolic Dehn fillings.
Moreover, we show that determining these finitely many manifolds
is decidable.  
\end{abstract}



\maketitle
%
%

\section{Introduction}
Thurston demonstrated that if one has a hyperbolic knot complement,
all but finitely many Dehn fillings give hyperbolic manifolds \cite{Th}. The 
example with the largest known number of non-hyperbolic Dehn 
fillings is the figure-eight knot complement, which has 10 
fillings which are not hyperbolic. It is conjectured that this
is the maximal number that can occur. 
Call a Dehn filling {\it exceptional} if it is not hyperbolic.
Previous authors have distinguished between hyperbolike manifolds, which are irreducible, atoroidal, with infinite fundamental group,
and hyperbolic manifolds, defining an exceptional filling to be one
which is not hyperbolike.  But by the geometrization theorem \cite{Per02, Per03,KleinerLott06,  MorganTian07, CaoZhu06}, a manifold is hyperbolic if and only if it is hyperbolike, so we need no longer
make this distinction. 
Bleiler and Hodgson \cite{BH} showed that there are at most 24 exceptional 
Dehn fillings, using Gromov and Thurston's $2\pi$-theorem and
estimates on cusp size due to Colin Adams \cite{A}. We made an
improvement on the $2\pi$-theorem \cite{Ag1}, independently discovered by Lackenby \cite{L}, and used improved lower bounds
on cusp size due to Cao and Meyerhoff \cite{CM}, to get an upper bound
of 12 exceptional Dehn fillings. In this paper, we
show that there are at most finitely many one cusped hyperbolic
manifolds which have exceptional Dehn fillings $r_1, r_2$
such that $\D(r_1,r_2)>5$. This theorem is sharp, in that
the $(-2,3,8)$ pretzel link complement $W'$ admits two exceptional Dehn fillings
$r_1,r_2$ with $\D(r_1,r_2)=5$. By hyperbolic Dehn filling
the other cusp of $W'$, we see that there are infinitely many 
3-manifolds with $\D(r_1,r_2)=5$. This theorem implies
that there are only finitely many one cusped hyperbolic 3-manifolds
with $>8$ exceptional Dehn fillings, since there can be at most $8$ curves on a torus with $\D(r_1,r_2)\leq 5$.  We also 
prove the existence of an algorithm which will determine the manifolds which have $\D(r_1,r_2)>5$ for exceptional Dehn fillings $r_1,r_2$, and
therefore which manifolds may have $>8$ exceptional Dehn fillings. It is conjectured
that there are only four orientable hyperbolic manifolds with two
exceptional Dehn fillings $r_1,r_2$ so that $\D(r_1,r_2)>5$. 
These manifolds are the Dehn fillings $W(-1), W(5), W(5/2)$, and $W(-2)$
on one component of the Whitehead link complement $W$, and
are conjectured by Gordon to be the only such examples \cite[Conjecture 3.4]{G}. Of these
four examples, only
the figure eight knot complement $W(-1)$ has $>8$ exceptional
Dehn fillings. It is conjectured that there are only finitely many 
one cusped hyperbolic 3-manifolds with $>6$ exceptional Dehn
fillings \cite[Problem 1.77]{Kirby}. 

Our theorem depends on results of Anderson,
Canary, Culler and Shalen \cite{ACCS} which have been subsequently strengthened by a culmination of various
results in Kleinian groups. In particular, we use Ohshika's result that
Kleinian groups are limits of geometrically finite Kleinian groups \cite{ohshika-2005} (see also \cite{NamaziSouto}),
which depends on the classification of Kleinian groups, including
the tameness conjecture \cite{Ag04, CG06} (see also \cite{Soma06}) and the ending lamination conjecture \cite{BCM, Minsky03}, and 
generalizes many previous results on density of geometrically
finite Kleinian groups in the space of Kleinian groups \cite{Bromberg07, BrockBromberg04}, going back
to work of Jorgensen on the space of punctured torus groups \cite{Jorgensen03} and
Thurston's double limit theorem \cite{Thurstonfiber, Otal96}. We remark that
Ohshika's argument \cite{ohshika-2005} depends on the ending
lamination conjecture for general Kleinian groups, which has 
been claimed by Brock, Canary and Minsky and disseminated 
through talks at various conferences, but for which there is as
yet no preprint. The preprints \cite{BCM, Minsky03} only treat
the case of freely-indecomposable Kleinian groups, but we need
the result for freely decomposable groups for the application in this paper.  
The key geometric consequence that we make use of is the following result stated in Theorem \ref{cuspvol}:
if a maximal horocusp in an orientable hyperbolic manifold has volume
$< \pi-\e$ for $\e>0$, then the volume of the manifold is uniformly bounded as a function of $\e$. 

\section{Background and definitions}
Let $N$ be a hyperbolic 3-manifold. Then $N \cong \HH^3/\G$
where $\G\cong \pi_1(N)$ is a discrete torsion-free subgroup of $\Isom(\HH^3)$, with $p:\HH^3\to N$ the covering map.
We will assume that the reader is familiar with Margulis' constant $\e$
and the thick-thin decomposition in the case of hyperbolic 3-manifolds;
see \cite[Chapter 5]{Th} for an introduction to this. The following 
definitions are needed only for part of the proof of Lemma \ref{freebicuspid}.
If $g$ is a loxodromic isometry of hyperbolic 3-space $\HH^3$, we shall let $A_g$ denote the
hyperbolic geodesic which is the axis of $g$. The cylinder about $A_g$ of radius $r$ is the open
set $Z_r(g) = \{x \in \HH^3 |\  dist(x,A_g) < r\}$. If
$C$ is a simple closed geodesic in $N$ then there is a primitive loxodromic isometry $g \in \G$ with
$p(A_g/\langle g\rangle) = C$. For any $r > 0$,  the projection $p(Z_r(g)/\langle g\rangle)$ of $Z_r(g)$ under the covering projection
is a neighborhood of $C$ in $N$. For sufficiently small $r > 0$ we have
$\{h \in  \G |\ h(Z_r(g)) \cap Z_r(g) \neq \emptyset\} = \langle g\rangle$.
Let $R$ denote the supremum of the set of $r$ for which this condition holds. We define
$\tube(C) = Z_R(g)/\langle g\rangle$ to be the maximal tube about $C$. 

\begin{figure}[htb] 
	\begin{center}
	\psfrag{r}{$\HH^3$}
	\psfrag{h}{$\tilde{\mathcal{H}}$}
	\psfrag{j}{$\g(\tilde{\mathcal{H}})$}
	\psfrag{g}{$\g^{-1}(\tilde{\mathcal{H}})$}
	\psfrag{a}{$\a$}
	\psfrag{b}{$\b$}
	\psfrag{x}{$x$}
	\psfrag{y}{$\g^{-1}(x)$}
	\epsfig{file=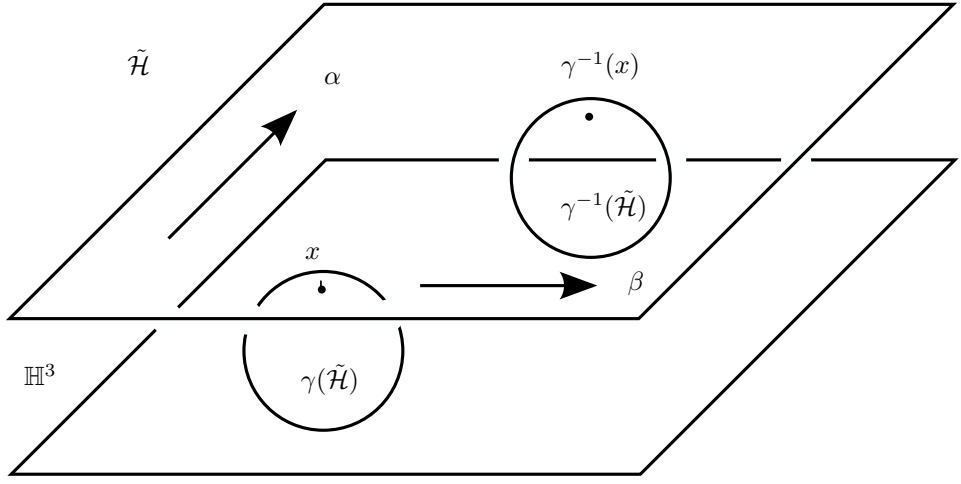, width=\textwidth}
	\caption{\label{bicuspid} A bicuspid group associated to a maximal horocusp}
	\end{center}
\end{figure}

There is a similar situation when we have a maximal parabolic subgroup 
$P< \G$. We define an open horoball to be a subset of 
$\HH^3$ isometric to $ \{ (z, t) | \ z\in \CC,  t>1\}$, where $\HH^3=\{(z,t)|\ t>0\} \subset \CC\times \RR$
is the upper half-space model of hyperbolic 3-space. 
There is an open horoball $\tilde{\mathcal{H}}\subset \HH^3$
so that $P= \{ g\in \G |\ g(\tilde{\mathcal{H}}) = \tilde{\mathcal{H}}\}$. 
Choose $\tilde{\mathcal{H}}$ maximal subject to the condition that
$\{ h\in \G |\ h(\tilde{\mathcal{H}})\cap \tilde{\mathcal{H}} \neq \emptyset\}= P$. 
Then we call $p(\tilde{\mathcal{H}})=\mathcal{H} \cong \tilde{\mathcal{H}}/P$ a {\it maximal horocusp}. If $P\cong \ZZ^2$, then $\mathcal{H}\cong S^1\times S^1\times \RR$, and if $P\cong \ZZ$, then $\mathcal{H}\cong S^1\times \RR^2$. Since $\mathcal{H}$ is maximal subject to 
this condition, there exists $\g\in \G-P$ so that 
$\g(\d \tilde{\mathcal{H}})\cap \d \tilde{\mathcal{H}} =x$. Note that
$\g^{-1}(\d \tilde{\mathcal{H}})\cap \d\tilde{\mathcal{H}}=\g^{-1}(x)$ will also
be a point of tangency. Then $\g^{-1}(\tilde{\mathcal{H}})$ is also a horoball
of height one call the Adam's horoball \cite{A} (see Figure \ref{bicuspid}).
We will call the group $B= \langle P, \g\rangle$ a {\it bicuspid} group,
so called because it corresponds to a point in $\HH^3$ where
the closure of two maximal horoballs meet, corresponding to two preimages
of a neighborhood of a cusp of $N$. In this paper, when
$\ZZ^2\cong P$, we will use the notation $P =\langle \a, \b\rangle$,
where $\a, \b$ are parabolic elements. Then $B$
is {\it free bicuspid} if $B\cong P\ast \langle \g\rangle \cong \ZZ^2\ast \ZZ$. For any cusp of $N$ with associated maximal horocusp
$\mathcal{H}$, there will be finitely many bicuspid subgroups
of $\G$ up to conjugacy in one-to-one correspondence with
the $\G$ orbits of points $x\in \HH^3$ so that 
$x=\g(\d \tilde{\mathcal{H}})\cap \d \tilde{\mathcal{H}} $
for some $\g \in \G-P$. One may think of each of these points as 
``self-tangencies'' of $\d \mathcal{H}$. 

\section{Maximal embedded cusps of hyperbolic manifolds}
The following geometric result is made possible by recent
advances in Kleinian groups, and generalizes \cite[Lemma 4.3]{ACS06}. 
\begin{lemma} \label{freebicuspid}
Let $\langle \a, \b, \g \rangle$ be a free bicuspid group,
with maximal horocusp $\mathcal{H}$. Then $\Vol(\mathcal{H})\geq \pi$.
\end{lemma}
\begin{proof}
Let $N=\HH^3/\langle \a,\b, \g\rangle$. Since $\mathcal{H}$ is maximal, there is a point of self-contact on
$\d \mathcal{H}$. 
Let $\ZZ+\ZZ\cong \pi_{1}(\mathcal{H}) = \langle \a, \b \rangle < \pi_{1}(N)$. Since the cusp $\mathcal{H}$ is maximal,
there is a component of the preimage of $\mathcal{H}$
in the universal cover $\HH^3$ fixed by covering translations $\langle \a, \b\rangle$ such that 
 $\g (\partial \tilde{\mathcal{H}}) \cap \partial\tilde{\mathcal{H}} \neq \emptyset$. Denote $Q = \langle \a,\b,\g \rangle = \pi_{1}(N)$. 
We have $\Area(\d \mathcal{H})=\frac12 \Vol(\mathcal{H})$. 
For sake of contradiction, we may assume that $\Area(\d \mathcal{H})\leq 2\pi$. 
Choose $\a$ to represent the shortest path in $\d \mathcal{H}$,
and $\b$ the next shortest path. By \cite{A}, the length of a
geodesic path representing $\a$ in $\d\mathcal{H}$ is $\geq 1$. 
We may normalize $\a, \b, \g$ as elements in $\SL_2 \CC$ up to conjugacy, such that $\tilde{\mathcal{H}}$ is the horoball centered at $\infty$ of height 1, and 

$$\a = \left(\begin{array}{cc}1 & a \\0 & 1\end{array}\right), \b = \left(\begin{array}{cc}1 & b \\0 & 1\end{array}\right), \g=\left(\begin{array}{cc}c & -1 \\1 & 0\end{array}\right),$$
where $1\leq |a| \leq |b| \leq 4\pi /\sqrt{3}$, and $|c|\leq |b|$. 

To prove the theorem, first we approximate $Q$ by geometrically finite groups
$Q_i$ such that $Q_i = \langle \a_i, \b_i, \g_i \rangle$ and such
that $\a_i \to \a, \b_i \to \b, \g_i \to \g$. Denote by $N_i=\HH^3/Q_i$
the hyperbolic manifold which is homeomorphic to the interior of a boundary
connect sum $(T^2\times I)\amalg (S^1\times D^2) \cong \overline{N}_i$. Denote $S_i$ to be the torus boundary component of $\overline{N}_i$. 
In fact, one may assume that
in the groups $Q_i$, the only parabolic elements are conjugate into $\langle \a_i, \b_i\rangle$. This is possible by the
fact that geometrically finite groups without rank one parabolics are dense in the space of all Kleinian
groups in the algebraic topology \cite[Theorem 1.1]{ohshika-2005}. We may assume
that $N_i$ contains a maximal  embedded horocusp neighborhood $\mathcal{H}_i$ 
such that $\Vol(\mathcal{H}_i)\to\Vol(\mathcal{H})$. 

Next, we proceed as in Lemma 4.3 of \cite{ACS06}. 
Suppose that $(P_j)$ is an infinite sequence of distinct hyperbolic manifolds obtained
by Dehn filling $N_i$ along $S_i$ using Thurston's hyperbolic Dehn surgery theorem for geometrically
finite manifolds (see \cite{BonahonOtal88} or \cite{Comar:1996}).
Then $\pi_1(P_j)$ is free two generator, and the manifolds
$(P_j)$ converge geometrically to $N_i$. Moreover, 
the core curve of the Dehn filling $P_j$ 
of $N_i$ is isotopic to a geodesic $C_j$ in $P_j$.  The length $L_j$ of $C_j$ tends to 0 as $j \to\infty$;
and the sequence of maximal tubes $(\tube(C_j)), j\geq1$ converges geometrically to $\mathcal{H}_i$. In
particular $\underset{j \to \infty}{\lim} \Vol(\tube(C_j)) = \Vol(\mathcal{H}_i)$.
According to \cite[Corollary 4.2]{ACS06}, $\log 3$ is a strong Margulis number for each of the hyperbolic
manifolds $P_j$. It therefore follows from \cite[Corollary 10.5]{ACCS} that $\Vol( \tube(C_j)) > V (L_j)$,
where $V$ is an explicitly defined function such that $\underset{x\to 0} {\lim}V (x) = \pi$. In particular, this
shows that $\Vol(\mathcal{H}_i) \geq \underset{j\to \infty}{\lim}V (L_j) \geq \pi$. Now, since $\mathcal{H}_i$ converges geometrically 
to $\mathcal{H}$, we conclude that $\Vol(\mathcal{H})\geq \pi$. 
\end{proof}

{\bf Remark:} The estimates of \cite{ACCS} depend on
a paradoxical decomposition argument for a 2-generator free
group acting on $\HH^3$ given in \cite{CS92}.
We believe that it should be possible to give a more direct 
argument for the previous lemma by analyzing an appropriate
generalization of the paradoxical decomposition for free bicuspid
groups. 

{\bf Example:} Consider the group $\G=\langle \a,\b,\g\rangle$, where
$$\a = \left(\begin{array}{cc}1 & 4 \\0 & 1\end{array}\right), \b = \left(\begin{array}{cc}1 & 1+i\sqrt{3} \\0 & 1\end{array}\right), \g=\left(\begin{array}{cc} 2  & -1 \\1 & 0\end{array}\right).$$
One may show that this group is free. Consider the geodesic
planes in the upper half space model of $\HH^3$ bounding the
circles $\{ |z- c| = 2 |\ c \in  \ZZ 2+ \ZZ(1+i\sqrt{3})\}$, and
cut out the open half-spaces disjoint from $\infty$ bounded by these planes to obtain a region $R\subset \HH^3$. These
bounding circles have two orbits under the group $\langle \a, \b\rangle$,
with representatives at centers $c=0, 2$, and thus $R/\langle \a,\b\rangle$ will be homeomorphic to $T^2\times \RR$ with two geodesic
disks in its boundary. Under the map $\g$,
the circle $|z|=2$ is sent to the circle $|z-2|=2$, and therefore
$\g$ maps the corresponding planes bounding these circles to themselves. Thus,
$N=\HH^3/\G$ is obtained from $R/\langle \a, \b\rangle$
by gluing the two geodesic planes in the the boundary using the
identification given by $\g$. This has the effect of adding a
handle onto $T^2\times \RR$, and therefore $N$ is homeomorphic to the interior of a compression body with fundamental group identified
with $\G\cong \langle \a,\b\rangle \ast \langle \g\rangle$. Thus, we see
that $\G$ is free bicuspid. If $\tilde{\mathcal{H}}$ is the horosphere
centered at $\infty$ of height one in $\HH^3$, then $\tilde{\mathcal{H}}/\langle \a, \b\rangle$ embeds as a cusp $\mathcal{H}\subset N$, since it is disjoint from
the geodesic planes given above. We have  $\Vol(\mathcal{H}) = 2\sqrt{3}= 3.46\ldots $. Thus, the bound given in Lemma \ref{freebicuspid}
is fairly close to optimal. 

\begin{figure}[htb]
\begin{center}
 \epsfbox{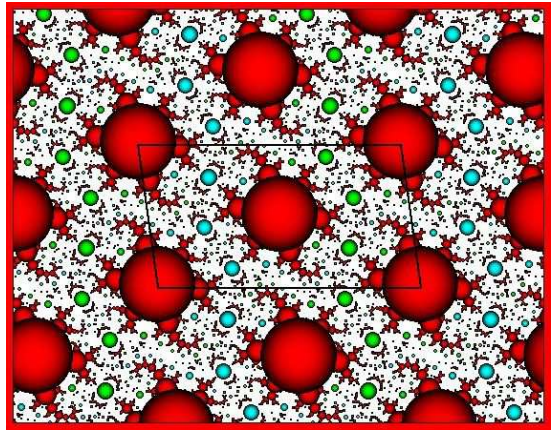}
 \caption{Pattern of horoballs in $\HH^3$ for the manifold v1902 3-fold cover 8 from the Snappea census with
cusp volume $3.238\ldots$
  }  \label{v1902}
\end{center}
\end{figure}

We did a search of  
manifolds in the Snappea census \cite{Snappea}, and found that the manifold
$M={\rm v1902}$ has an embedded cusp of volume $3.238\ldots$.
Moreover, the cusp group associated to {\rm v1902} is free.
One may see this by finding an irregular 3-fold cover $\tilde{M}\to M$ (cover
8 in Snappea's notation) which has 3 maximal cusps with
the same volume as $M$, and has $H_1(\tilde M)= \ZZ^3 + \ZZ/3\ZZ$ 
(see Figure \ref{v1902}). 
It follows that the bicuspid subgroup of $\pi_1 M$ is free,
since it lifts to $\pi_1 \tilde{M}$, and any bicuspid subgroup
of $\pi_1 \tilde{M}$ is free since ${\rm rank} H_1(\tilde{M};\ZZ/3\ZZ)=4$ (see 
the argument of \cite[Prop. 5.3]{ACS06}).

In fact, there will be an optimal constant
$\pi \leq R < 3.238\ldots$ such that the quotient manifold of every 
free bicuspid group will have an embedded horocusp $\mathcal{H}$ 
with $\Vol(\mathcal{H})\geq R$. It is an interesting question to 
obtain better upper and lower bounds on $R$. The following theorem 
would hold with $R$ replacing $\pi$.  

\begin{theorem} \label{cuspvol}
For $\e>0$, there exists a constant $V(\e)$ such that if $N$ is an orientable hyperbolic
3-manifold with a maximal horocusp $\mathcal{H}$ with $\Vol(\mathcal{H})< \pi-\e$,
then $\Vol(N)< V(\e)$. 
\end{theorem}
\begin{proof}
Let $Q< \pi_1 N$ be a bicuspid subgroup corresponding to the maximal horocusp $\mathcal{H}$. 
First, suppose that $Q$ is a free product $Q = \langle \a_1,\a_2\rangle \ast \langle \b\rangle \cong (\ZZ+\ZZ)\ast \ZZ$. 
Then $\Vol(\mathcal{H})>\pi $ by Lemma \ref{freebicuspid}.

Next, suppose $Q$ is not a free product. Then $Q$ must be
indecomposable, and in fact $\HH^3/Q$ must be finite volume. 
We get a contradiction in this case by taking a geometric limit. 
Suppose we have $\G_j< \PSL(2,\CC)$ with $N_j=\HH^3/\G_j$, such that $\Vol(N_j)\to \infty$ and maximal
bicuspid subgroups $Q_j < \G_j$, such that $Q_j$ is not free
bicuspid. Then $\Vol(\HH^3/Q_j)< \infty$ by \cite[Theorem IV.4.1]{JS79}, so $[\G_j:Q_j]<\infty$. We may assume that $\Vol(\HH^3/Q_j) \to \infty$ as $j\to \infty$,
since otherwise there exists $V$ such that $V>\Vol(\HH^3/Q_j)\geq \Vol(\HH^3/\G_j)$, which 
is what we wanted to conclude.  Normalize $Q_j$ as a standard
bicuspid group with $\mathcal{H}_j\subset \HH^3/Q_j$ the maximal cusp, and with self-tangency point of $x_j \in \d \mathcal{H}_j \subset  \HH^3/Q_j$.
Take a subsequence $j_k$ such that $(\HH^3/Q_{j_k}, x_{j_k})$ converges geometrically
to $(M_{\infty}, x_{\infty})$, where $\mathcal{H}_{\infty}\subset M_{\infty}$ is 
a maximal cusp with self-tangency point $x_{\infty}$ and $\Vol(M_{\infty})=\infty$. Then $Q_j$
converges algebraically to a bicuspid subgroup $Q_{\infty}< \pi_1(M_{\infty})$ such that
$\mathcal{H}_{\infty}$ lifts to a maximal cusp in $\HH^3/Q_{\infty}$.
By \cite[Theorem IV.4.1]{JS79}, $Q_{\infty} \cong (\ZZ+\ZZ)\ast \ZZ$,
so that $Q_{\infty}$ is a free bicuspid group. The cusps $\mathcal{H}_j \to \mathcal{H}_{\infty}$ in the
Gromov-Hausdorff topology. Since $\Vol(\mathcal{H}_{\infty} )\geq\pi$ by Lemma \ref{freebicuspid}, we see
that there is $k$ such that $\Vol(\mathcal{H}_{j_k} )> \pi-\e$. This gives a contradiction,
so we see that the assumption that $\Vol(\HH^3/Q_j) \to \infty$ is false, and thus there
exists $V(\e)$ such that $\Vol(N_j)<V(\e)$. 
\end{proof}

{\bf Remark:} A similar argument to the previous theorem shows
that there is a constant $V$ such that if $M$ is a hyperbolic 3-manifold
with Margulis constant $<\log 3$, then $\Vol(M)< V$ by applying
\cite{CS92}. It would be
interesting to get some idea of the distribution of Margulis constants $<\log 3$.

\section{Exceptional Dehn fillings}

In this section, we prove the main theorem of the paper. 
\begin{theorem} \label{distance}
There are only finitely many one cusped orientable
hyperbolic 3-manifolds of finite volume $N$ such that $N$ has
two exceptional Dehn fillings $a_1, a_2$ so that $\D(a_1,a_2)>5$. 
\end{theorem}
\begin{proof}
Fix some $\e < \pi-3$. If $\Vol(N)> V(\e)$, where 
$V(\e)$ is the constant from Theorem \ref{cuspvol}, 
then $N$ has a maximal horocusp $\mathcal{H}\subset N$
so that $\Vol(\mathcal{H})\geq \pi-\e > 3$. Then $\Area(\d\overline{\mathcal{H}}) > 6$. By the proof of \cite[Theorem 8.1]{Ag1}, 
if $a_1, a_2$ are two slopes on $\d \overline{N}$ so that
$\overline{N}(a_i)$ are exceptional Dehn fillings,
then $\D(a_1,a_2) \leq 6^2/\Area(\d \overline{\mathcal{H}}) <6$.
Thus, $N$ does not violate the theorem. 

Suppose the theorem is false. Then there is an infinite sequence
of orientable one cusped hyperbolic 3-manifolds $(N_i)_{i\in \NN}$ such that $i\to \infty$ and $N_i$ has two exceptional Dehn fillings
of distance $>5$. By the previous paragraph, we may assume that $\Vol(N_i) < V(\e)$.  
Let $\mathcal{H}_i$ be a maximal cusp neighborhood in $N_i$. 
By \cite[Theorem 5.12.1]{Th}, we may choose
a subsequence $J\subset \NN$ and an orientable hyperbolic 3-manifold
$M$ with $\geq 2$ rank two cusps so that each $N_{i}, i\in J$ is
obtained by hyperbolic Dehn filling on $M$, and therefore
$N_i$ converges to $M$ in the Gromov-Hausdorff topology. Moreover, $M$ has
a distinguished cusp with maximal cusp neighborhood $\mathcal{H}$
such that as $i\to \infty, i \in J$, $\mathcal{H}_i\to \mathcal{H}$
in the Gromov-Hausdorff topology (in fact, since $\Vol(\mathcal{H})\leq \pi-\e$,
we may assume that $M$ has at most three cusps). By \cite[Theorem 1.3]{Gordon98},
the distance between exceptional  filling slopes on the cusp
$\mathcal{H}$ of $M$
is $\leq 5$. Since $N_i \to M$ as $i\to \infty, i\in J$, it follows
that there is a subset $J_0\subset J$, $|J-J_0|<\infty$, so that the distance between
filling slopes   on $N_i$ for $i\in J_0$ is $\leq 5$.  
Let $\d \overline{M} = T_0 \cup T_1 \cup T_2$, where
$T_0$ corresponds to $\mathcal{H}$. 
For slopes $ q_j \in \QQ\cup \infty\cup \ast$, $j=0, 1, 2$,
let $M(q_0, q_1, q_2)$ be Dehn filling on the 
boundary component $T_i$ with slope $q_i$,  where
the $\ast$ denotes that the boundary component is 
unfilled. 
Let $\overline{N_i}= M(\ast, q_{i,1}, q_{i,2})$, $i\in J$. 
For each $0\leq j\leq 2$, there is a finite subset  $E_j\subset \QQ\cup \infty$
so that if $q_j\notin E_j$, then $M(q_0, q_1, q_2)$ is 
hyperbolic \cite[Theorem 5.8.2]{Th}. 
For $i\in J$ large enough, $q_{i,j}\notin E_j$, for $j=1,2$.
This implies that the only non-hyperbolic Dehn fillings
on $N_i$ must  correspond to a subset of the slopes $E_0$.
Let $E_h\subset E_0$ be such that  $t\in E_h$ if and only if 
$M(t,\ast,\ast)$ is hyperbolic. For $t_1, t_2 \in E_0-E_h$,
we have $\D(t_1,t_2)\leq 5$ by \cite[Theorem 1.3]{Gordon98} as noted above.
For each slope $t\in E_h$ there are subsets $ E_{t,j}\subset \QQ\cup \infty$
so that $M(t,q_1, q_2)$ is hyperbolic if $q_j \notin E_{t,j}$ (again
by \cite[Theorem 5.8.2]{Th}).
Let $F_j = E_j \cup_{t\in E_h} E_{t,j}$. Then if
$q_{i,j} \notin F_j$ and $t\in E_h$,
then $N_i(t)$ is also hyperbolic. Therefore, we see that
for $i$ large, the distance between two exceptional slopes of $N_i$ is at
most 5. This gives a contradiction to our assumption
that there is a subsequence of manifolds contradicting
the theorem with volumes bounded by $V(\e)$, thus
proving the theorem. 
\end{proof}

\begin{corollary} \label{finite}
There are only finitely many one cusped orientable
hyperbolic 3-manifolds of finite volume $N$ such that $N$ has
$>8$ exceptional Dehn fillings. 
\end{corollary}
\begin{proof}
By \cite[Lemma 8.2]{Ag1}, if the distance between exceptional
Dehn fillings is $\leq 5$, take the next largest prime and add
one, to conclude there are at most
$8$ exceptional Dehn fillings on $N$. 
\end{proof}

\section{An algorithm to find exceptional Dehn fillings}
We would like to classify the finitely many manifolds
with two exceptional Dehn fillings of distance $>5$  given by Theorem \ref{distance}.
Since the proof of Theorem \ref{cuspvol} is by contradiction,
it's not clear that there is a procedure one could run which would
for a given $\e$ identify all of the manifolds which have 
an embedded cusp of volume $<\pi-\e$. Thus, we must devise a method
to classify manifolds with a small volume cusp. This is 
similar to the procedures implemented in \cite{GMT, GMM07},
which find manifolds with small radius tubes around a short
geodesic, or small volume cusped manifolds. We modify their approach
to show that finding the manifolds with exceptional fillings in
Theorem \ref{distance} and Corollary \ref{finite} is decidable. 
First, we need some preliminary results. 

\begin{theorem} \label{exceptional}
Given a finite volume hyperbolic manifold $M$, there is
an algorithm which will determine the set of all exceptional Dehn
fillings on $M$ for which any proper sub-Dehn filling is hyperbolic. 
\end{theorem}
\begin{proof}
We approximate the hyperbolic structure
on the cusped manifold well enough to determine
the rough shape (for example using interval arithmetic ) of simultaneously embedded equal sized  horocusp neighborhoods
of all the cusps. This uses the algorithms to compute Ford domains 
in \cite{Ri83, HW89, Weeks93} (see also \cite{Manning}). We then apply \cite{Ag1, L} to determine which slopes on each cusp have length $\leq 6$. We perform all of the Dehn
fillings along these short slopes, determining which ones result in hyperbolic manifolds
using the algorithms described in \cite{Ri83, Weeks93, Manning},
then repeat. At each stage we get hyperbolic manifolds with
fewer cusps, and so this process eventually terminates
with a finite collection of cusped hyperbolic manifolds along
with a finite collection of slopes associated to each cusp, so that
every exceptional Dehn filling on $M$ is obtained
by Dehn filling on one of the slopes associated to one of the members of this finite collection. 
\end{proof}

\begin{lemma} \label{isoperimetric}
Let $N$ be a hyperbolic 3-manifold, and $R\subset N$
be a connected region, such that $im\{\pi_1(R)\to \pi_1(N)\}$ is
elementary. Then $\Vol(R)\leq \frac12 \Area(\d R)$. 
\end{lemma}
\begin{proof}
This follows from the isoperimetric inequality for $\HH^3$
plus the fact that the image of $\pi_1 R$ in $\pi_1 N$ is 
amenable. 
Alternatively, one may also use the 
method of \cite[Lemma 3.2, Theorem 4.1]{Rafalski07}, which is 
essentially a calibration argument.
\end{proof}

\begin{lemma} \label{relatorbound}
Suppose $Q=\langle \a, \b, \g\rangle$ is a (discrete torsion-free) bicuspid group which is not free. Let $N=\HH^3/Q$.  Let $w(x, y, z) \in (\langle x\rangle \times \langle y\rangle) \ast \langle z\rangle$ be a cyclically reduced word, 
such that $w(\a, \b, \g)=1 \in Q$. Let
$d(w)$ be the number of occurrences of $z^{\pm 1}$ in the word $w$.
Then $\Vol(\HH^3/Q) \leq \pi( d(w)-2)$. 
\end{lemma}
\begin{proof}
This result generalizes  \cite{Cooper99}.
Consider the word $w= r_1 z^{\e_1} r_2 z^{\e_2}\cdots r_k z^{\e_j}$,
where we assume that $\e_i \neq 0$, and $r_i \in \langle x\rangle \times \langle y\rangle -\{1\}$. Then
$d(w) = \sum_{i=1}^j |\e_i|$. For $k=0, \ldots, d(w)$, let 
$w_k(x,y,z)$ be a prefix of $w$ so that $d(w_k)=k$ and $w_k$
ends in $z^{\pm 1}$ (thus $w_{d(w)}=w$). Then 
we may find a string of geodesics $\g_k$, $k=1,\ldots, d(w)$ in $\HH^3$ so 
that $\g_k$ connects $w_{k-1}(\a,\b,\g) (\infty)$ and $w_k(\a,\b,\g)(\infty)$, and so that $\g_{d(w)}$ connects $w_{d(w)-1}(\a,\b,\g)(\infty)$ to $\infty=w(\infty)$. We may find a map of a disk $r: D\to \HH^3$ with $d(w)$
punctures in $\d D$ so that $r(\d D) \subset \g_1 \cup \cdots \cup \g_{d(w)}$ by coning each $\g_i$ to $\infty$, for $1< i<d(w)$.
The disk $r(D)$ will be made of $d(w)-2$ triangles, and therefore $\Area(r(D))=\pi(d(w)-2)$.
 Under the projection map $p: \HH^3 \to N$, each $\g_i$ will project
to the same geodesic $\g\subset N$, and the boundary projects to 
$p\circ r(\d D) \subset \g$.

Consider a component $U \subset N \backslash p\circ r(D)$, and let $W\subset U$ be a compact 
submanifold with boundary which is a deformation retract of $U$. For
each sphere $S\subset N\backslash W$, there is a ball
$B\subset N, \d B=S$. If $W\subset B$ for some ball $B$, then $im\{ \pi_1 W\to \pi_1 N\}=im\{\pi_1 U \to \pi_1 N\}=1$. Otherwise, $B\subset N\backslash W$
for each 2-sphere $S\subset N\backslash W$, and therefore $N\backslash W$
is irreducible. If $T^2\subset N\backslash W$ is a torus, then
either $T^2$ is incompressible in $N$, and therefore $T=\d V$,
where $V$ is a neighborhood of a cusp end of $N$, or $T$
is compressible in $N$. If $T$ is compressible in $N$, then $T=\d V$,
where $V$ is a ball with knotted hole (possibly a solid torus).
In either case, we see that  if $W\subset V$, then $im\{ \pi_1 U\to \pi_1 N\}$ is contained in an elementary subgroup. So if $im\{\pi_1 U\to \pi_1 N\}$ is not contained in an elementary subgroup, we 
conclude that $N\backslash W$ is irreducible and atoroidal. 
This implies that some component $\S\subset \d W$ must
have $\chi(\S)<0$, since any sphere or torus boundary component
of $\d W$ would bound a ball or ball with knotted hole component or cusp neighborhood
of $N\backslash W $, which would contradict the fact that $N\backslash U\supset p(r(D))$
is connected. Therefore, $\chi(N\backslash W)<0$. 
 Let
$G=im\{ Q \to \pi_1(N\backslash W)\}$, where we identify 
$Q$ with the image of $\pi_1 (\g \cup \mathcal{H})$, where $\mathcal{H}\subset N\backslash W$ is a horocusp. We then have $G$ must
be a free bicuspid subgroup of $\pi_1 (N\backslash W)$,
since $\chi(N \backslash W)<0$ by \cite[Theorem IV.4.1]{JS79}. But by hypothesis,
$w(x, y, z)$ gives a non-trivial relation for $G$, since
$p(r(D))\subset N\backslash W$. This gives a contradiction,
unless for each component $U$ of $N\backslash r(p(D))$,
we have $im\{ \pi_1 U \to \pi_1 N\}$ lies in an elementary
subgroup. Therefore, we conclude that $\Vol(U)\leq\frac12 \Area(\d U)$ by Lemma \ref{isoperimetric}. We conclude that $\Vol(N) = \Vol(N\backslash r(p(D))) \leq \Area(r(p(D)))\leq \pi (d(w)-2)$.
\end{proof}

\begin{theorem} \label{volumelist}
Given a rational number $V$, there is an algorithm which will
find a finite collection $M_1, \ldots, M_m$ of finite volume
orientable hyperbolic 3-manifolds so that any hyperbolic
3-manifold of volume $<V$ is obtained by 
Dehn filling on one of the manifolds $M_i$.
\end{theorem}
\begin{proof}
Let $\varepsilon$ be Margulis' constant for hyperbolic 3-manifolds.
By the method of proof of the Jorgensen-Thurston theorem,
there is a constant $C$ so that if $M$ is a hyperbolic 3-manifold
with $\Vol(M)<V$, then $M_{thick(\varepsilon)}$ admits
a triangulation with $<CV$ tetrahedra. The first step of
the algorithm is to take $< CV$ tetrahedra, and glue
them together in all possible ways to get an orientable
manifold with Euler characteristic zero. Next, we run the algorithm
described in \cite{HW89, Weeks93, Manning} to decide which of these manifolds has a hyperbolic interior. 
Then all hyperbolic 3-manifolds of volume $<V$
will be obtained by Dehn filling on one of the resulting finite
collection of manifolds. 
\end{proof}

\begin{theorem} \label{algorithm}
There is an algorithm which will determine the finitely many one cusped orientable hyperbolic 3-manifolds
$N_1, \ldots, N_k$ such that $N_i$ has two exceptional Dehn
fillings of distance $>5$. 
\end{theorem}
\begin{proof}
We know that if $N$ is a one-cusped orientable hyperbolic
3-manifold which has an embedded maximal horocusp neighborhood $\mathcal{H}\subset N$ such that $\Vol(\mathcal{H})>3$,
then any two Dehn fillings on $N$ have distance $\leq 5$.
Thus, we must determine the one-cusped orientable 
hyperbolic 3-manifolds $N$ which have a maximal cusp of volume $\leq 3$.
Such manifolds have a non-free bicuspid subgroup $Q=\langle \a, \b,\g\rangle$ of finite index in $\pi_1 N = \G$, such that $\mathcal{H}$
lifts to a maximal cusp of $\HH^3/Q$. Thus, we must determine
the bicuspid groups $Q$ which have a maximal cusp of $\HH^3/Q$
of volume $\leq 3$. By Theorem \ref{cuspvol}, we know that
$\Vol(\HH^3/Q)<  V(\e)$, for $\e<\pi-3$. But we don't know how to compute
$V(\e)$ explicitly. So we first show that there is an algorithm
which will determine all bicuspid groups $Q$ such that
$\HH^3/Q$ has a maximal cusp of volume $<\pi-\e$. 

As in the proof of Lemma \ref{freebicuspid}, we normalize
a general bicuspid group $Q=\langle \a, \b, \g\rangle$ so that
$$\a = \left(\begin{array}{cc}1 & a \\0 & 1\end{array}\right), \b = \left(\begin{array}{cc}1 & b \\0 & 1\end{array}\right), \g=\left(\begin{array}{cc}c & -1 \\1 & 0\end{array}\right),$$
where $1\leq |a| \leq |b| \leq 2A /\sqrt{3}$, and $|c|\leq |b|$,
where $A=\Area(\d \mathcal{H})$. This gives a finite parameter
space $\mathcal{P}$ to search through, and we know that when $A<2\pi$,
the group $Q$ must either be non-free bicuspid, or else
$Q$ is indiscrete. Furthermore, if $Q$ is discrete but
not a free product, then either $Q$ has finite covolume,
or $Q$ has torsion. In this case, there will be a non-trivial
word $w(x,y,z) \in (\langle x \rangle \times \langle y \rangle) \ast \langle z\rangle \cong (\ZZ+\ZZ)\ast \ZZ$ 
such that $w(\a, \b, \g)=1 \in Q$. A lift of $w(\a,\b,\g)$ to
$\SL_2\CC$ gives a matrix 
$W(a,b,c)=\left(\begin{array}{cc}\ast  & \ast  \\ p(a,b,c) & \ast \end{array}\right),$ such that $W(a,b,c)=\left(\begin{array}{cc} 1&0\\ 0&1\end{array}\right)$ when the variables $(a,b,c)$ correspond to the group $Q$. 
The lower left entry of $W(a,b,c)$ may be regarded as a polynomial
$p(a,b,c)\in \ZZ[a,b,c]$ which vanishes on the parameters corresponding
to $Q$. 

If 
$\mu = \left(\begin{array}{cc}w  & x \\ y & z\end{array}\right)$
represents a matrix in $\SL_2 \CC$, and if $\tilde{\mathcal{H}}$
is the horoball of height $1$ centered at $\infty$ in $\HH^3$,
then $\mu(\tilde{\mathcal{H}})$ is a horoball of height 
$1/|y|^2$ for $y\neq 0$. Thus, if $\tilde{\mathcal{H}}\cap \mu(\tilde{\mathcal{H}})=\emptyset$, we see that $|y|\geq 1$. So, if $\mu \in Q$
is an element of the bicuspid group $Q$, and $\tilde{\mathcal{H}}$ is precisely invariant under $Q$, then we see that
if $|y| <1$, we must have $\tilde{\mathcal{H}} = \mu(\tilde{\mathcal{H}}) $, and therefore $\mu$ is a parabolic element, so $y=0$.

Similarly, if $Q$ is indiscrete, then the closure $\overline{Q}\leq \SL_2 \CC$
must be dense in a Lie subgroup of $\SL_2\CC$, and therefore
there are elements $w(\a,\b, \g)\in Q$ such that $|W(a,b,c)-I|^2<\e$,
for any small $\e$. In particular, if 
$$W(a,b,c) = \left(\begin{array}{cc}\ast & \ast \\ p(a,b,c) & \ast \end{array}\right),$$
where $p(a,b,c) < 1$, then either $p(a,b,c)=0$, or the three generator
group associated to the parameters $(a,b,c)$ does not have
an embedded horocusp $\mathcal{H}\subset \HH^3/\langle \a, \b,\g\rangle$ which is a projection of a horoball at height 1 in $\HH^3$, and is therefore not bicuspid with respect to these generators. 
Find a finite collection of such polynomials 
$\{ p_1, \ldots, p_k\} \in \ZZ[a,b,c]$ so that the sets
$\{ (a,b,c) | p_i(a,b,c)< 1 \}$ cover $\mathcal{P}$. 
Each such polynomial is determined by a word 
$w_i(x,y,z)$. We may find a finite collection since $\mathcal{P}$
is compact. To make this search algorithmic, enumerate
$p(a,b,c)$ for all $w(x,y,z) \in \ZZ^2\ast \ZZ$. 
Cover $\mathcal{P}$ by compact subsets,
such as cubes with dyadic vertices and sidelength $2^{-n}$. For a given polynomial $p(a,b,c)$,
and a given cube $C$, one may determine whether  
$|p(a,b,c)|<1$ for all points $(a,b,c)\in C$, by determining
the maximal value of $p(a,b,c)$ on $C$, which is algorithmic
(an interior maximum is determined by computing a point
where the gradient $\nabla p(a,b,c)=0$, which may be
computed via algebraic geometry, whereas boundary
maxima may be determined by Lagrange multipliers 
inductively on the faces of the cube). Inductively, we 
refine the covering and increase the number of polynomials,
until we find a covering of $\mathcal{P}$ (this is similar to 
the process employed in \cite[Proposition 1.28]{GMT}).

Points $(a,b,c)\in\mathcal{P}$ where 
$p_i(a,b,c)=0$ correspond to irreducible bicuspid representations of 
$\ZZ^2\ast \ZZ$. If $0<p_i(a,b,c)< 1$, then either the group
is indiscrete, or it does not have a maximal cusp $\mathcal{H}$
normalized as above. Thus, all of the discrete non-free
bicuspid groups with a cusp of volume $\leq 3$ will occur
somewhere in the parameter space $\mathcal{P}$ with some 
$w_i(\a,\b,\g)=1$, and therefore $p_i(a,b,c)=0$. 
By Lemma \ref{relatorbound}, the
covolume of the discrete torsion-free $Q$ will be bounded by the complexity $\pi (d(w_i)-2)$. 
Let $V$ be the supremum of these numbers.
Thus, we have shown the existence of an algorithm
to determine a bound on the covolume of bicuspid
groups with a maximal cusp of volume $\leq 3$. 

The next step is to enumerate all orientable hyperbolic
manifolds with volume $< V$. We use Theorem \ref{volumelist}
to give a finite collection $M_1,\ldots, M_m$ of finite
volume hyperbolic 3-manifolds so that any hyperbolic
3-manifold of volume $<V$ is obtained by hyperbolic
Dehn filling on one of the manifolds $M_j$. 
Finally, we use Theorem \ref{exceptional} to determine
the set of all exceptional Dehn fillings on each $M_j$.
We then search this list of exceptional Dehn fillings
for any one cusped manifolds with two exceptional
fillings of distance $>5$. By Theorem \ref{distance},
there will be only finitely many isometry types of 
hyperbolic one cusped manifolds with two exceptional
Dehn fillings of distance $>5$. 
\end{proof}

\section{Conclusion}
\begin{figure}[htb]
  \subfigure{\epsfig{figure=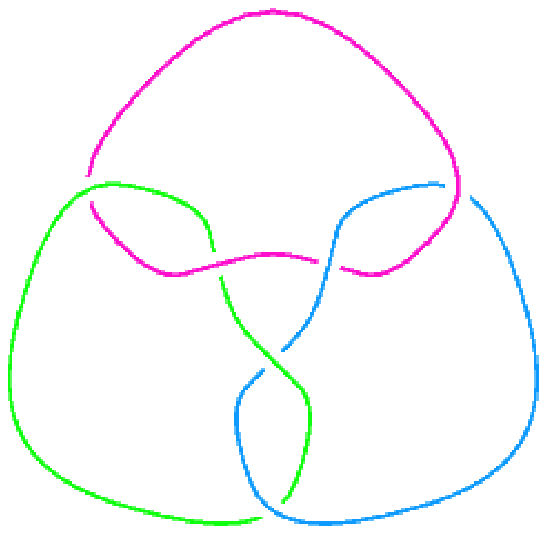,angle=0,width=.45\textwidth}}\quad
    \subfigure{\epsfig{figure=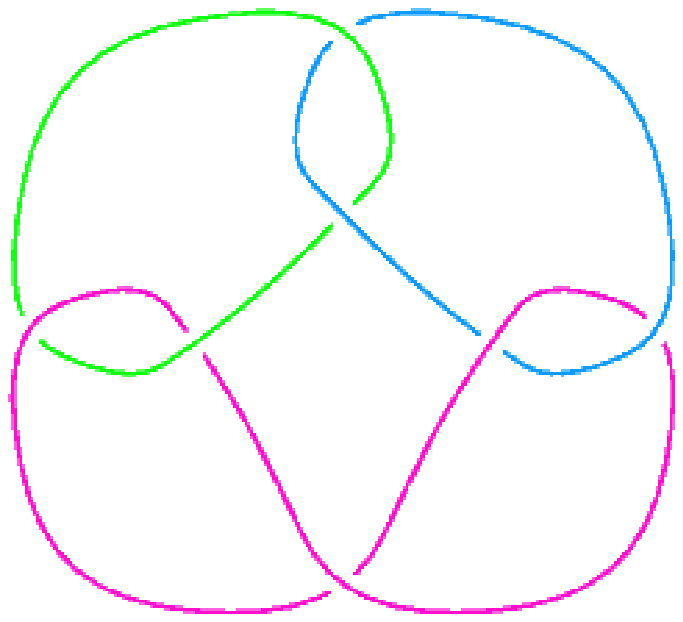,angle=0,width=.45\textwidth}}\\
\subfigure{\epsfig{figure=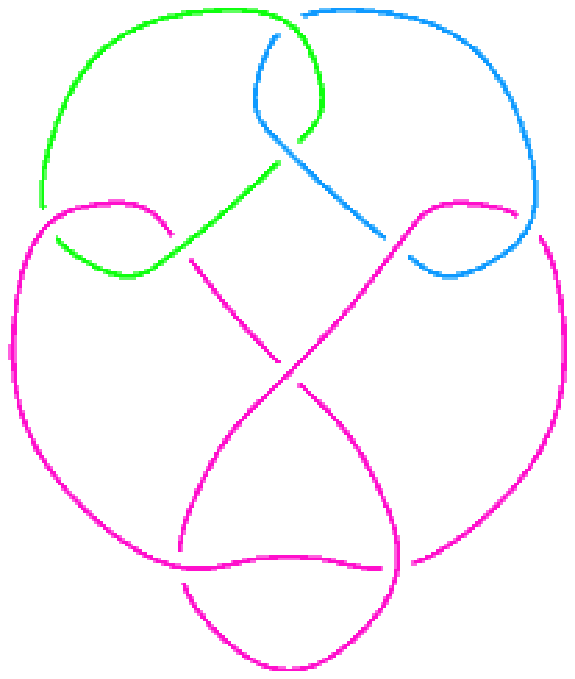,angle=0,width=.45\textwidth}}
\quad 
\subfigure{\epsfig{figure=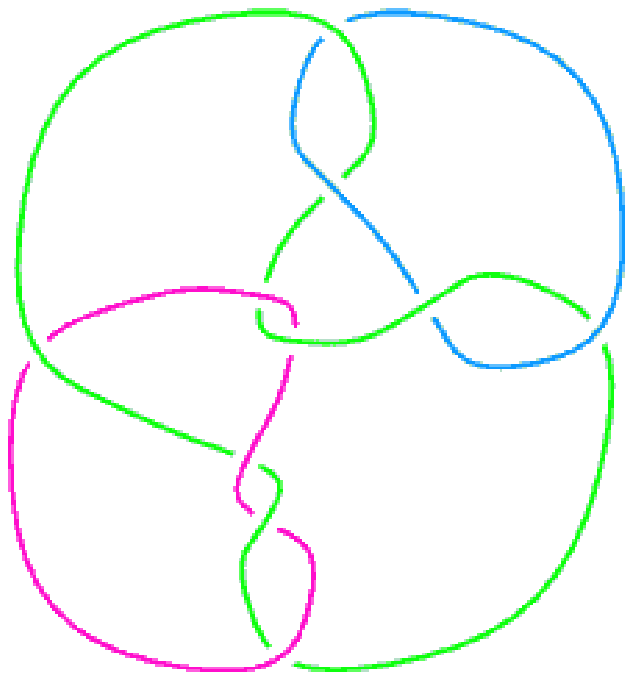,angle=0,width=.45\textwidth}}
    \caption{\label{smallcusps} Some 3 component links with cusps of volume $\leq 3$}
\end{figure}

There are several problems suggested by the results in this
paper. Of course, there is Gordon's conjectured classification of one cusped manifolds with two exceptional fillings of distance $>5$ 
(see \cite[Problem 1.77, Conjecture (B)]{Kirby}). 
As we point out, it will suffice to classify the manifolds with
a cusp of volume $\leq 3$. Searching through the Snappea
census of cusped orientable manifolds \cite{Snappea}, there appears to
be few one or two-cusped manifolds which have cusp volume $ \leq 3$ and
which are not obtained by Dehn filling on a manifold
with corresponding cusp of volume very close to $3$ and which have
fewer than three cusps. The collection of 3-cusped manifolds
with a cusp of volume $\leq 3$ is finite. 
We show several examples of three cusped
manifolds with a cusp of volume $\leq 3$ in Figure \ref{smallcusps},
and which give rise to infinite families of 3-manifolds with
maximal cusps of volume $\leq 3$ via Dehn filling.
We remark that if one wants to classify manifolds with $>8$
exceptional Dehn fillings, it suffices to classify manifolds with
two exceptional Dehn fillings $r_1, r_2$ with $\D(r_1,r_2)\leq 6$.
Thus, it suffices to consider manifolds with a cusp of volume
$\leq 2\frac47$. This would significantly simplify the search space in  the algorithm
proposed in Theorem \ref{algorithm}.

A related question is to understand the structure of non-free
bicuspid groups. One may obtain examples by considering
genus 2 manifolds in which one of the compression bodies
$C$ such that $\d_+ C= \S_2$, $\d_- C=T^2$. It is natural
to ask whether this construction gives all examples. We
conjecture that all hyperbolic
manifolds with bicuspid fundamental group are genus 2. If so, this would help the search 
algorithm in Theorem \ref{algorithm} since one would
need only add relators coming from embedded curves
on the boundary of a genus two compression body to
obtain all of the bicuspid groups with small cusp volume.

It would be nice to go beyond Gordon's conjecture, and to 
classify all finite-volume orientable hyperbolic manifolds with
two exceptional Dehn fillings of distance $>4$. We suspect
that the only examples are obtained from Dehn filling on the
$(-2,3,8)$ pretzel link and the Whitehead link. To achieve this using the methods in 
this paper, one would have to classify hyperbolic manifolds
with cusp volume $<4.5$, and this does not seem feasible.
To make further progress using geometric methods, one will
likely have to delve further into the geometry of cusped
manifolds, and try to understand Ford domains of cusped
manifolds with small volume horocusps. For each slope in a cusp, one may consider
the minimal area of a complete surface inside the Ford domain which is asymptotic to  the given slope. If this area is $>2\pi$, then the Dehn
filling along this slope cannot be exceptional. It may be
possible to numerically estimate these areas for a fixed
type of Ford domain, and could possibly give better 
estimates on the number of exceptional surgeries than
the estimates obtained from cusp volume. 
It seems likely that the methods of \cite{GMM07}
should be suited to this sort of analysis, and it may be possible
to resolve \cite[Problem 1.77]{Kirby}.

%
%
\bibliographystyle{../hamsplain}
\bibliography{../refs}


%
%
\end{document}